\documentclass[12pt,twoside]{article}
\usepackage{amssymb}
\usepackage{amsmath}
\begin{document}
\setcounter{page}{1}
\newtheorem{t1}{Theorem}[section]
\newtheorem{d1}{Definition}[section]
\newtheorem{c1}{Corollary}[section]
\newtheorem{l1}{Lemma}[section]
\newtheorem{r1}{Remark}[section]

\newcommand{\cA}{{\cal A}}
\newcommand{\cB}{{\cal B}}
\newcommand{\cC}{{\cal C}}
\newcommand{\cD}{{\cal D}}
\newcommand{\cE}{{\cal E}}
\newcommand{\cF}{{\cal F}}
\newcommand{\cG}{{\cal G}}
\newcommand{\cH}{{\cal H}}
\newcommand{\cI}{{\cal I}}
\newcommand{\cJ}{{\cal J}}
\newcommand{\cK}{{\cal K}}
\newcommand{\cL}{{\cal L}}
\newcommand{\cM}{{\cal M}}
\newcommand{\cN}{{\cal N}}
\newcommand{\cO}{{\cal O}}
\newcommand{\cP}{{\cal P}}
\newcommand{\cQ}{{\cal Q}}
\newcommand{\cR}{{\cal R}}
\newcommand{\cS}{{\cal S}}
\newcommand{\cT}{{\cal T}}
\newcommand{\cU}{{\cal U}}
\newcommand{\cV}{{\cal V}}
\newcommand{\cX}{{\cal X}}
\newcommand{\cW}{{\cal W}}
\newcommand{\cY}{{\cal Y}}
\newcommand{\cZ}{{\cal Z}}

\def\cl{\centerline}
\def\bd{\begin{description}}
\def\be{\begin{enumerate}}
\def\ben{\begin{equation}}
\def\benn{\begin{equation*}}
\def\een{\end{equation}}
\def\eenn{\end{equation*}}
\def\benr{\begin{eqnarray}}
\def\eenr{\end{eqnarray}}
\def\benrr{\begin{eqnarray*}}
\def\eenrr{\end{eqnarray*}}
\def\ed{\end{description}}
\def\ee{\end{enumerate}}

\def\al{\alpha}
\def\b{\beta}
\def\bR{\bar\R}
\def\bc{\begin{center}}
\def\ec{\end{center}}
\def\d{\dot}
\def\D{\Delta}
\def\del{\delta}
\def\ep{\epsilon}
\def\g{\gamma}
\def\G{\Gamma}
\def\h{\hat}
\def\iny{\infty}
\def\La{\Longrightarrow}
\def\la{\lambda}
\def\m{\mu}
\def\n{\nu}
\def\noi{\noindent}
\def\Om{\Omega}
\def\om{\omega}
\def\p{\psi}
\def\pr{\prime}
\def\r{\ref}
\def\R{{\bf R}}
\def\ra{\rightarrow}
\def\s{\sum_{i=1}^n}
\def\si{\sigma}
\def\Si{\Sigma}
\def\t{\tau}
\def\th{\theta}
\def\Th{\Theta}
\def\vep{\varepsilon}
\def\vp{\varphi}
\def\pa{\partial}
\def\un{\underline}
\def\ov{\overline}
\def\fr{\frac}
\def\sq{\sqrt}
\def\WW{\begin{stack}{\circle \\ W}\end{stack}}
\def\ww{\begin{stack}{\circle \\ w}\end{stack}}
\def\st{\stackrel}
\def\Ra{\Rightarrow}
\def\bi{\begin{itemize}}
\def\ei{\end{itemize}}
\def\i{\item}
\def\bt{\begin{tabular}}
\def\et{\end{tabular}}
\def\lf{\leftarrow}
\def\nn{\nonumber}
\def\va{\vartheta}
\def\wh{\widehat}
\def\vs{\vspace}
\def\Lam{\Lambda}
\def\sm{\setminus}
\def\su{\sum^\iny_{n=1}}
\def\ba{\begin{array}}
\def\ea{\end{array}}
{\baselineskip 15truept

\bc
{\bf  Some applications of the generalized Bernardi - Libera - Livingston integral operator on univalent functions}\\
\ec

\bc
{\bf  M. Eshaghi Gordji} \\
Department of Mathematics, Faculty of Science\\
Semnan University, Semnan, Iran  \\
e-mail: madjid.eshaghi@gmail.com \\

{\bf A. Ebadian}\\
Department of Mathematics, Faculty of Science \\
Urmia University, Urmia, Iran\\
e-mail: a.ebadian@mail.urmia.ac.ir\\

\ec
}

\baselineskip 24truept

\begin{abstract}
In this paper by making use of the generalized Bernardi - Libera - Livingston integral operator
 we introduce and study some new subclasses of univalent functions.  Also we investigate the
 relations between those classes and the classes which are studied by Jin-Lin Liu.
\end{abstract}

\noi {\bf Key Words} : {\it Starlike, convex, close-to-convex, quasi-convex, strongly starlike, strongly convex functions}.

\noi 2005 Mathematics Subject Classification : Primary 30C45,
30C75.

\section{Introduction}

Let $A$ be the class of functions of the form, $ f(z) = z + \sum^\iny_{n=2} a_nz^n$
which are analytic in the unit disk $U = \{ z : |z| < 1 \}$, also let $S$ denote the
subclass of $A$ consisting of all univalent functions in $U$.  Suppose $\la$ is a real
number with $0 \le \la < 1$, $A$ function $f \in S$ is said to be starlike of
order $\la$ if and only if $Re \left\{ \fr{zf'(z)}{f(z)} \right\} > \la, z \in U$,
also $f \in S$ is said to be convex of order $\la$ if and only if $Re \left\{ 1 + \fr{zf''(z)}{f'(z)} \right\} > \la, z \in U$,
we denote by $S^*(\la), C(\la)$ the classes of starlike and convex functions of order $\la$ respectively.
It is well known that $f \in C(\la)$ if and only if $zf' \in S^*(\la)$.  Let $f \in A$ and $g \in S^*(\la)$
then $f \in K(\b, \la)$ if and only if $Re \left\{ \fr{zf'(z)}{g(z)} \right\} > \b, z \in U$ where $0 \le \b < 1$.
 These functions are called close-to-convex functions of order $\b$ type $\la$.  A function $f \in A$ is called quasi-convex of
 order $\b$ type $\la$ if there exists a function $g \in C(\la)$ such that $Re \left\{ \fr{(zf'(z))'}{g'(z)} \right\} > \b$.
  We denote this class by $K^*(\b, \la)$ [10].  It is easy to see that $f \in K^*(\b, \g)$ if and only if $zf' \in K(\b, \g)$ [9].  For $f \in A$ if for some $\la (0 \le \la < 1)$ and $\eta(0 < \eta \le 1)$ we have
$$ \left| arg \left( \fr{zf'(z)}{f(z)} - \la \right) \right| < \fr{\pi}{2} \eta , ~~ z \in U \eqno{(1.1)}$$
then $f(z)$ is said to be strongly starlike of order $\eta$ and type $\la$ in $U$ and we denote this class by $S^*(\eta, \la)$.  If $f \in A$ satisfies the condition
$$ \left| arg \left( 1 + \fr{zf''(z)}{f'(z)} - \la \right) \right| < \fr{\pi}{2} \eta , ~~ z \in U \eqno{(1.2)}$$
for some $\la$ and $\eta$ as above then we say that $f(z)$ is strongly convex of order $\eta$ and type $\la$ in $U$ and we denote this class by $C(\eta, \la)$.
Clearly $f \in C(\eta, \la)$ if and only if $zf' \in S^*(\eta, \la)$, specially we have $S^*(1, \la) = S^*(\la)$ and $C(1, \la) = C(\la)$.

For $c > -1$ and $f \in A$ the generalized Bernardi - Libera - Livingston integral operator $L_cf$ is defined as follows
$$L_cf(z) = \fr{c+1}{z^c} ~~ \int\limits^z_0 t^{c-1} f(t)dt. \eqno{(1.3)}$$
This operator for $c \in N = \{ 1, 2, 3, \cdots \}$ was studied by
Bernardi [1] and for $c = 1$ by Libera [5] (see also [8]).  Now by
making use of the operator given by (1.3) we introduce the following
classes. \benrr
S^*_c(\la) & = & \{ f \in A : L_cf \in S^*(\la)\} \\
C_c(\la) & = & \{ f \in A : L_cf \in C(\la)\} \\
K_c(\b, \la) & = & \{ f \in A : L_cf \in K(\b, \la)\} \\
K^*_c(\b, \la) & = & \{ f \in A : L_cf \in K^*(\b, \la)\} \\
ST_c(\eta, \la) & = & \left\{ f \in A : L_cf \in S^*(\eta, \la), \fr{z(L_cf(z))'}{L_cf(z)} \ne \la, z \in U \right\} \\
CV_c(\eta, \la) & = & \left\{ f \in A : L_cf \in C(\eta, \la),
\fr{(z(L_cf(z))')'}{(L_cf(z))'} \ne \la, z \in U \right\}. \eenrr
Obviously $f \in CV_c(\eta, \la)$ if and only if $zf' \in ST_c(\eta,
\la)$.  J. L. Liu [6] and [7] introduced and investigated similarly
the classes $S^*_\si(\la), C_\si(\la), K_\si(\b, \la), K^*_\si(\b,
\la), ST_\si(\eta, \la), CV_\si(\eta, \la)$  by making use of the
integral operator $I^\si f$ given by
$$I^\si f(z) = \fr{2^\si}{z\G(\si)} \int\limits^z_0 \left( \log \fr{z}{t} \right)^{\si -1} f(t)dt, \si > 0, f \in A. \eqno{(1.4)}$$
The operator $I^\si$ is introduced by Jung, Kim and Srivastava [3]
and then investigated by Uralogaddi and Somanatha [12],  Li [4] and
Liu [6].  For the integral operators given by (1.3) and (1.4) we
have  easily verified following relationships.
$$I^\si f(z) = z + \sum^\iny_{n=2} \left( \fr{2}{n+1} \right)^\si a_nz^n \eqno{(1.5)}$$
$$L_cf(z) = z + \sum^\iny_{n=2} \fr{c+1}{n+c} a_nz^n \eqno{(1.6)}$$
$$z(I^\si L_cf(z))' = (c + 1) I^\si f(z) - cI^\si L_c f(z) \eqno{(1.7)}$$
$$z(L_cI^\si f(z))' = (c+1) I^\si f(z) - c L_c I^\si f(z) . \eqno{(1.8)} $$
It follows from (1.5) that one can define the operator $I^\si$ for
any real number $\si$. In this paper we investigate the properties
of the classes $S^*_c(\la), C_c(\la), K_c(\b, \la), K^*_c(\b,
\la),\\ ST_c(\eta, \la), CV_c(\eta, \la)$, also we study the
relations between these classes by the classes which are introduced
by Liu in [6] and [7].  For our purposes we need the following
lemmas.

\noi {\bf Lemma 1.1} [9].  Let $u = u_1 + iu_2, v = v_1 + iv_2$ and
let $\psi(u, v)$ be a complex function $\psi : D \subset \mathbb{C}
\times \mathbb{C} \ra \mathbb{C}$.  Suppose that $\psi$ satisfies
the following conditions

(i) $\psi(u, v) $ is continuous in $D$

(ii) $(1, 0) \in D$ and $Re \{ \psi(1, 0)\} > 0$

(iii) $Re \{ \psi(iu_2, v_1)\} \le 0$ for all $(iu_2, v_1) \in D$ with $v_1 \le - \fr{1 + u^2_2}{2}$.

Let $p(z) = 1 + \sum\limits^\iny_{n=2} c_nz^n$ be analytic in $U$ so that $(p(z), zp'(z)) \in D$ for all $z \in U$.  If $Re \{ \psi (p(z), zp'(z))\} > 0, z \in U$ then $Re \{ p(z) \} > 0, z \in U. $

\noi {\bf Lemma 1.2} [11].  Let the function $p(z) = 1 +
\sum\limits^\iny_{n=1} c_nz^n$ be analytic in $U$ and $p(z) \ne 0, z
\in U$ if there exists a point $z_0 \in U$ such that $|arg (p(z))| <
\fr{\pi}{2} \eta$ for $|z| < |z_0|$ and $|arg~ p(z_0)| = \fr{\pi}{2}
\eta$ where $0 < \eta \le 1$ then $\fr{z_0 p'(z_0)}{p(z_0)} =
ik\eta$ and $k \ge \fr{1}{2} (r + \fr{1}{r})$ when $arg~ p(z_0) =
\fr{\pi}{2} \eta$ also $k \le \fr{-1}{2} (r + \fr{1}{r})$ when $arg~
p(z_0) = \fr{-\pi}{2} \eta$, and $p(z_0)^{1/\eta} = \pm ir (r > 0)$.

\section{Main Results}

In this section we obtain some inclusion theorems.

\noi {\bf Theorem 2.1} : (i) For $f \in A$ if $Re \left\{
\fr{zf'(z)}{f(z)} - \fr{z(L_cf(z))'}{L_cf(z)} \right\} > 0$, then
$S^*_c(\la) \subset S^*_{c+1}(\la)$.

(ii) For $f \in A$ if $Re \left\{ \fr{zf'(z)}{f(z)} - \fr{z(L_{c+1} f(z))'}{L_{c+1} f(z)} \right\} > 0$ then $S^*_{c+1}(\la) \subset S^*_c(\la)$.

\noi {\bf Proof} : (i) Suppose that $f \in S^*_c(\la)$ and set
$$ \fr{z(L_{c+1} f(z))'}{L_{c+1} f(z)} - \la = (1 - \la) p(z) \eqno{(2.1)}$$
where $p(z) = 1 + \sum\limits^\iny_{n=2} c_nz^n$.  An easy calculation shows that
$$ \fr{ \fr{z(L_{c+1} f(z))'}{L_{c+1} f(z)} \left[ 2 + c + \fr{z(L_{c+1} f(z))''}{(L_{c+1} f(z))'} \right]}{ \fr{z(L_{c+1} f(z))'}{L_{c+1} f(z)} + c + 1 } = \fr{zf'(z)}{f(z)} .\eqno{(2.2)}$$
By setting $H(z) = \fr{z(L_{c+1} f(z))'}{L_{c+1} f(z)}$ we have
$$ 1 + \fr{z(L_{c+1} f(z))''}{(L_{c+1} f(z))'} = H(z) + \fr{zH'(z)}{H(z)}. \eqno{(2.3)}$$
By making use of (2.3) in (2.2) since $H(z) = \la + (1 - \la) p(z)$ so we obtain
$$ (1 - \la) p(z) + \fr{(1 - \la) zp'(z)}{\la + c + 1 + (1 - \la) p(z)} = \fr{zf'(z)}{f(z)} - \la. \eqno{(2.4)}$$
If we consider $\psi(u, v) = (1 - \la)u + \fr{(1 - \la)v}{\la + c + 1 + (1 - \la)u}$ then $\psi(u, v)$ is a continuous function in $D = \left\{ \mathbb{C} - \fr{\la + c + 1}{\la - 1} \right\} \times \mathbb{C}$ and $(1, 0) \in D$ also $\psi(1, 0) > 0$ and for all $(iu_2, v_1) \in D$ with $v_1 \le - \fr{1 + u^2_2}{2}$ we have
$$Re~~ \psi(iu_2, v_1) = \fr{(1 - \la)(\la + c + 1)v_1}{(1-\la)^2u^2_2+(\la+c+1)^2} \le \fr{-(1-\la)(\la+c+1)(1 + u^2_2)}{2[(1-\la)^2u^2_2 + (\la+c+1)^2]} < 0.$$
Therefore the function $\psi(u, v)$ satisfies the conditions of
Lemma 1.1 and since in view of the assumption by considering
(2.4) we have $Re\{ \psi(p(z), zp'(z))\} > 0$ therefore Lemma 1
implies that $Re~ p(z) > 0, z \in U$ and this completes the proof.

(ii) For  proving this part of theorem by the same method and using the easily verified formula similar to (2.2) by replacing $c+1$ with $c$ we get the desired result.

\noi {\bf Theorem 2.2} : (i) For $f \in A$ if $Re \left\{
\fr{zf'(z)}{f(z)} - \fr{z(L_cf(z))'}{L_cf(z)} \right\} > 0$ then
$C_c(\la) \subset C_{c+1}(\la)$.

(ii) For $f \in A$ if $Re \left\{ \fr{zf'(z)}{f(z)} - \fr{z(L_{c+1} f(z))'}{L_{c+1}f(z)} \right\} > 0$ then $C_{c+1}(\la) \subset C_c(\la)$.

\noi {\bf Proof} : (i) In view of part (i) of Theorem 1 we can write\\ $f \in C_c(\la) \Leftrightarrow L_cf \in C(\la) \Leftrightarrow z(L_cf)' \in S^*(\la) \Leftrightarrow L_c zf' \in S^*(\la) \Leftrightarrow zf' \in S^*_c(\la) \Ra zf' \in S^*_{c+1} (\la) \Leftrightarrow L_{c+1} zf' \in S^*(\la) \Leftrightarrow z(L_{c+1} f)' \in S^*(\la) \Leftrightarrow L_{c+1} f \in C(\la) \Leftrightarrow f \in C_{c+1}(\la). $

By the similar way we can prove the part (ii) of theorem.

\noi {\bf Theorem 2.3} : If $c \ge - \la$ then $f \in S^*(\la)$
implies $f \in S^*_c(\la)$.

\noi {\bf Proof} :  By differentiating logarithmically from both sides of (1.3) with respect to $z$ we obtain
$$ \fr{z(L_c f(z))'}{L_cf(z)} + c = \fr{(c+1) f(z)}{L_c f(z)} . \eqno{(2.5)}$$
Again differentiating logarithmically from both sides of (2.5) we have
$$ p(z) + \fr{zp'(z)}{c + \la + p(z)} = \fr{zf'(z)}{f(z)} - \la \eqno{(2.6)}$$
where $p(z) = \fr{z(L_cf(z))'}{L_cf(z)} - \la$.  Let us consider $\psi(u, v) = u + \fr{v}{u + c + \la}$, then $\psi$ is a continuous function in $D = \{ \mathbb{C} - (- c - \la)\} \times \mathbb{C}$ and $(1, 0) \in D$ also $Re ~ \psi(1, 0) > 0$.  If $(iu_2, v_1) \in D$ with $v_1 \le - \fr{1 + u^2_2}{2}$ then $Re ~ \psi(iu_2, v_1) = \fr{v_1(c + \la)}{u^2_2 + (c + \la)^2} \le 0$, also since $f \in S^*(\la)$ then (2.6) gives $Re(\psi(p(z), zp'(z))) = Re \left\{ \fr{zf'(z)}{f(z)} - \la \right\} > 0$.  Therefore Lemma 1 concludes that $Re\{p(z)\}>0$ and this completes the proof.

\noi {\bf Corollary 2.4} :  If $c \ge \la$ then $f \in C(\la)$
implies $f \in C_c(\la)$.

\noi {\bf Proof} : We have \\$f \in C(\la) \Leftrightarrow zf' \in S^*(\la) \La zf' \in S^*_c(\la) \Leftrightarrow L_c zf' \in S^*(\la) \Leftrightarrow z(L_cf)' \in S^*(\la) \Leftrightarrow L_cf \in C(\la) \Leftrightarrow f \in C_c(\la)$.

\noi {\bf Theorem 2.5} : (i) For $f \in A$ if $\left| arg \left(
\fr{zf'(z)}{f(z)} - \la \right) \right| \le \left| arg \left(
\fr{z(L_c f(z))'}{L_cf(z)} - \la \right) \right|, z \in U$ then
$ST_c(\eta, \la) \subset ST_{c+1} (\eta, \la), c > -1$.

(ii) For $f \in A$ if  $\left| arg \left( \fr{zf'(z)}{f(z)} - \la \right) \right| \le \left| arg \left( \fr{z(L_{c+1} f(z))'}{L_{c+1}f(z)} - \la \right) \right|, z \in U$ then $ST_{c+1} (\eta, \la) \subset ST_c(\eta, \la), c > -1$.

\noi {\bf Proof} : (i) Let $f \in ST_c(\eta, \la)$ and put
$$ \fr{z(L_{c+1} f(z))'}{L_{c+1} f(z)} = \la + (1 - \la) p(z) \eqno{(2.7)}$$
where $p(z) = 1 + \sum\limits^\iny_{n=1} c_nz^n$ is analytic in $U$ with $p(z) \ne 0, z \in U$.  It is easy to see that
$$ z(L_{c+1} f(z))' + (c+1) L_{c+1} f(z) = (c + 2) f(z) . \eqno{(2.8)} $$
Differentiating logarithmically with respect to $z$ from both sides of (2.8) gives
$$ \fr{z \left( \fr{z(L_{c+1} f(z))'}{L_{c+1} f(z)} \right)'}{\fr{z(L_{c+1} f(z))'}{L_{c+1} f(z)} + c + 1} + \fr{z(L_{c+1} f(z))'}{L_{c+1} f(z)} = \fr{zf'(z)}{f(z)} . \eqno{(2.9)} $$
Now by making use of (2.7) in (2.9) we have
$$ \fr{(1 - \la) zp'(z)}{\la + c + 1 + (1 - \la) p(z)} + (1 - \la) p(z) = \fr{zf'(z)}{f(z)} - \la. \eqno{(2.10)}$$
Suppose that there exists $z_0 \in U$ in such a way $|arg(p(z))|
< \fr{\pi}{2} \eta$ for $|z| < |z_0|$ and $|arg(p(z_0))| =
\fr{\pi}{2} \eta$, then by Lemma 1.2 we have
$\fr{z_0p'(z_0)}{p(z_0)} = ik\eta$ and $p(z_0)^{1/\eta} = \pm
ir(r > 0)$ where $k \ge \fr{1}{2} (r + \fr{1}{r}) $ when
$arg(p(z_0)) = \fr{\pi}{2} \eta$ and $k \le \fr{-1}{2} (r +
\fr{1}{r})$ when $arg(p(z_0)) = \fr{-\pi}{2} \eta$.  If
$p(z_0)^{1/\eta} = ir$ then $arg(p(z_0)) = \fr{\pi}{2} \eta$ and
by considering (2.10) we have \benrr
& & \left| arg \left( \fr{z_0(L_c f(z_0))'}{L_cf(z_0)} - \la \right)\right| \ge arg \left( \fr{z_0f'(z_0)}{f(z_0)} - \la \right) \\
& & = arg \left\{ (1 - \la) p(z_0) \left[ 1 + \fr{ik\eta}{\la + c + 1 + (1-\la) r^\eta e^{i \fr{\pi}{2} \eta}} \right] \right\} \\
& & = \fr{\pi}{2} \eta \\
& & + \tan^{-1} \left\{ \fr{k \eta [\la + c + 1 + r^\eta (1 - \la) \cos \fr{\pi}{2} \eta ]}{(\la + c + 1)^2 + r^{2\eta}(1-\la)^2 + (1-\la) (\la+c+1) \cos \fr{\pi}{2} \eta + k \eta r^\eta(1-\la) \sin \fr{\pi}{2} \eta} \right\} \\
& & \ge \fr{\pi}{2} \eta ~~ (\mbox{Because}~~ k \ge \fr{1}{2} (r + \fr{1}{r}) \ge 1)
\eenrr
which is a contradiction by $f(z) \in ST_c(\eta, \la)$.

Now suppose $p(z_0)^{1/\eta} = - ir$ then $arg(p(z_0)) = \fr{-\pi}{2} \eta$ and we have
\benrr
& &  - \left| arg \left( \fr{z_0(L_c f(z_0))'}{L_cf(z_0)} - \la \right)\right| \le arg \left( \fr{z_0f'(z_0)}{f(z_0)} - \la \right) \\
& & = \fr{-\pi}{2} \eta +  arg \left\{ 1 +  \fr{ik\eta}{\la + c + 1 + (1-\la) r^\eta e^{-i \fr{\pi}{2} \eta}}  \right\} \\
& & = \fr{-\pi}{2} \eta \\
& & + \tan^{-1} \left\{ \fr{k \eta [\la + c + 1 + r^\eta (1 - \la) \cos \fr{\pi}{2} \eta ]}{(\la + c + 1)^2 + r^{2\eta}(1-\la)^2 + 2r^\eta(1-\la) (\la+c+1) \cos \fr{\pi}{2} \eta - k \eta r^\eta(1-\la) \sin \fr{\pi}{2} \eta} \right\} \\
& & \le \fr{-\pi}{2} \eta ~~ (\mbox{Because}~~ k \le \fr{-1}{2} (r + \fr{1}{r}) \le -1)
\eenrr
which contradicts our assumption that $f \in ST_c(\eta, \la)$, therefore $|arg(p(z))| < \fr{\pi}{2}, z \in U$ and finally $\left| arg \left( \fr{z(L_{c+1} f(z))'}{L_{c+1} f(z)} - \la \right) \right| < \fr{\pi}{2} \eta, z \in U$.  However since for every $\la (0 \le \la < 1)$ we have $\fr{z(L_{c+1} f(z))'}{L_{c+1} f(z)} \ne \la$ thus we have $f \in ST_{c+1} (\eta, \la)$ and the proof is complete.

(ii) The proof of this part of theorem is similar with the proof of part (i) and we omit the proof.

\noi {\bf Corollary 2.6} : (i) For $f \in A$ if $\left|arg \left(
\fr{zf'(z)}{f(z)} - \la \right) \right| \le \left| arg \left(
\fr{z(L_c f(z))'}{L_cf(z)} - \la \right) \right|, z \in U$ then
$CV_c (\eta, \la) \subset CV_{c+1} (\eta, \la)$.

(ii) For $f \in A$ if, $\left| arg \left( \fr{zf'(z)}{f(z)} - \la \right) \right| \le \left| arg \left( \fr{z(L_{c+1} f(z))'}{L_{c+1} f(z)} - \la \right) \right| , z \in U$ then we have \\ $CV_{c+1}(\eta, \la) \subset CV_c(\eta, \la)$.

\noi {\bf Proof} : We give only the proof of part (i) and for this we have \\$f \in CV_c(\eta, \la) \Leftrightarrow L_c f \in C(\eta, \la) \Leftrightarrow z(L_cf)' \in S^*(\eta, \la) \Leftrightarrow L_c zf' \in S^*(\eta, \la) \Leftrightarrow zf' \in ST_c(\eta, \la) \La zf' \in ST_{c+1}(\eta, \la) \Leftrightarrow L_{c+1} zf' \in S^*(\eta, \la) \Leftrightarrow z(L_{c+1}f)' \in S^*(\eta, \la) \Leftrightarrow L_{c+1} f \in C(\eta, \la) \Leftrightarrow f \in CV_{c+1}(\eta, \la)$.

\noi {\bf Theorem 2.7} : For every $c > -1$ we have $CV_c(\eta,
\la) \subset ST_c(\eta, \la)$.

\noi {\bf Proof} : Let $f \in CV_c(\eta, \la)$ then $\left| arg \left( 1 + \fr{z(L_cf(z))''}{(L_c f(z))'} - \la \right) \right| < \fr{\pi}{2} \eta, z \in U$ and $\fr{(z(L_c f(z))')'}{(L_cf(z))'} \ne \la, z \in U$.  Suppose that
$$ \fr{z(L_c f(z))'}{L_cf(z)} = \la + (1 - \la) p(z) \eqno{(2.11)}$$
where $p(z) = 1 + \sum\limits^\iny_{n=2} c_nz^n$ is analytic in $U$ with $p(z) \ne 0$ for all $z \in U$.  Differentiating both sides of (2.11) logarithmically with respect to $z$ gives
$$ 1 + \fr{z(L_c f(z))''}{(L_c f(z))'} - \la = (1 - \la) p(z) + \fr{(1 - \la) zp'(z)}{\la + (1 - \la) p(z)}. $$
If there exists a point $z_0 \in U$ such that  $|arg(p(z))| < \fr{\pi}{2} \eta (|z| < |z_0|)$ and $|arg (p(z_0))| = \fr{\pi}{2} \eta$ then by Lemma 2 we obtain $\fr{z_0p'(z_0)}{p(z_0)} = ik\eta$ and $p(z_0)^{1/\eta} = \pm ir(r > 0)$ where $k \ge \fr{1}{2} (r + \fr{1}{r})$  when $arg(p(z_0)) = \fr{\pi}{2} \eta$ and $k \le - \fr{1}{2} (r + \fr{1}{r})$ when $arg(p(z_0)) = - \fr{\pi}{2} \eta$.  Suppose that $arg(p(z_0)) = \fr{-\pi}{2} \eta$ then
\benrr
& & arg \left\{ 1 + \fr{z_0(L_cf(z_0))''}{(L_c f(z_0))'} - \la \right\} \\
& & = arg \left\{ (1 - \la)r^\eta e^{-i \fr{\pi}{2}\eta} \left[ 1 + \fr{ik\eta}{\la + (1-\la)r^\eta e^{-i \fr{\pi}{2} \eta}} \right] \right\} \\
& & = \fr{-\pi}{2} \eta + arg \left\{ 1 + \fr{ik\eta}{\la + (1 - \la) r^\eta e^{-i \fr{\pi}{2} \eta}} \right\} \\
& & = \fr{-\pi}{2} \eta + \tan^{-1} \left\{ \fr{k\eta [\la + (1 - \la) r^\eta \cos \fr{\pi}{2} \eta]}{\la^2 + 2\la(1-\la)r^\eta \cos \fr{\pi}{2} \eta + (1-\la)^2 r^{2\eta} - k\eta(1 - \la)r^\eta \sin \fr{\pi}{2} \eta} \right\} \\
& & \le \fr{-\pi}{2} \eta ~~~ (\mbox{Because}~~ k \le \fr{-1}{2} (r + \fr{1}{r}) \le - 1)
\eenrr
which is a contradiction by $f \in CV_c(\eta, \la)$.  For the case $arg(p(z_0)) = \fr{\pi}{2}\eta$ by the same way and considering $k \ge \fr{1}{2} (r + \fr{1}{r}) \ge 1$ we obtain
$$arg \left\{ 1 + \fr{z_0(L_cf(z_0))''}{(L_cf(z_0))'} - \la \right\} \ge - \fr{\pi}{2} \eta.$$
This also contradicts our assumption that $f \in CV_c(\eta, \la)$, thus we have $|arg(p(z))| < \fr{\pi}{2} \eta (z \in U)$ and finally
$$ \left| arg \left( \fr{z(L_cf(z))'}{L_cf(z)} - \la \right) \right| < \fr{\pi}{2} \eta, ~~ z \in U.$$

\noi {\bf Theorem 2.8} : (i) If for every $f \in A$ and $g \in
S^*_c(\la)$ we have
$$Re \left\{ \fr{z \fr{d}{dz} \left( \fr{L_c zf'(z)}{L_cg(z)}\right)}{\fr{z(L_cg(z))'}{L_cg(z)} + c} \right\} > 0 \eqno{(2.12)} $$
and
$$ Re \left\{ \fr{zg'(z)}{g(z)}- \fr{z(L_cg(z))'}{L_cg(z)} \right\} > 0 \eqno{(2.13)} $$
then $K_c(\b, \la) \subset K_{c+1}(\b, \la)$.

(ii) If for every $f \in A$ and $g \in S^*(\la)$ we have
$$Re \left\{ \fr{z \fr{d}{dz} \left( \fr{L_{c+1} zf'(z)}{L_{c+1}g(z)}\right)}{\fr{z(L_{c+1}g(z))'}{L_{c+1}g(z)} + c} \right\} > 0 \eqno{(2.14)} $$
and
$$ Re \left\{ \fr{zg'(z)}{g(z)}- \fr{z(L_{c+1}g(z))'}{L_{c+1}g(z)} \right\} > 0 \eqno{(2.15)} $$
then $K_{c+1}(\b, \la) \subset K_{c}(\b, \la)$.

\noi {\bf Proof} : (i) Let $f \in K_c (\b, \la)$ then there exists a function $\vp(z) \in S^*(\la)$ such that
$$Re \left \{ \fr{z(L_c f(z))'}{\vp(z)} \right\} > \b,~~ z \in U.$$
There is a function $g$ in such a way $L_cg(z) = \vp(z)$ therefore $g \in S^*_c(\la)$ and we have $Re \left\{ \fr{z(L_cf(z))'}{L_cg(z)} \right\} > \b, z \in U$.  Suppose that
$$ \fr{z(L_{c+1} f(z))'}{L_{c+1} g(z)} - \b = (1 - \b) p(z) \eqno{(2.16)} $$
where $p(z) = 1 + \sum\limits^\iny_{n=1} c_nz^n$.  Now in view of (2.12) we can write
\benrr
& & 0 > Re \left\{ \fr{- z \fr{d}{dz} \left( \fr{L_c zf'(z)}{L_cg(z)} \right)}{\fr{z(L_cg(z))'}{L_cg(z)} + c }\right\} = Re \left\{ \fr{zL_c zf'(z) (L_cg(z))' - z(L_c zf'(z)' L_cg(z)}{L_cg(z)[z(L_cg(z))' + cL_c g(z)} \right\} \\
& & = Re \left\{ \fr{z(L_cf(z))' [z(L_cg(z))' + cL_cg(z)] - L_cg(z) [z(L_c zf'(z))' + cL_c z f'(z)]}{L_cg(z) [z(L_c g(z))' + cL_cg(z)]} \right\} \\
& & = Re \left\{ \fr{z(L_cf(z))'}{L_cg(z)} \right\} - Re \left\{ \fr{cL_czf'(z) + z(L_czf'(z))'}{z(L_cg(z))' + cL_cg(z))} \right\} .
\eenrr
Therefore we have
$$ Re \left\{ \fr{z(L_cf(z))'}{L_cg(z)} \right\} < Re \left\{ \fr{z(L_czf'(z))' + c(L_czf'(z))}{z(L_cg(z))' + cL_cg(z)} \right\} \eqno{(2.17)} $$
Now by easy computation we obtain the following identities.
$$z(L_c zf'(z))' + c(L_czf'(z)) = \fr{c+1}{c+2} [z(L_{c+1} zf'(z))' + (c+1) (L_{c+1} zf'(z))] \eqno{(2.18)} $$
$$z(L_c g(z))' + c(L_cg(z) = \fr{c+1}{c+2} [z(L_{c+1} g(z))' + (c+1) (L_{c+1} g(z))] .\eqno{(2.19)} $$
By making use of (2.18) and (2.19) in (2.17) we get
\benrr
& & Re \left\{ \fr{z(L_c f(z))'}{L_cg(z)} \right\} < Re \fr{z(L_{c+1} zf'(z))' + (c+1) (L_{c+1} zf'(z))}{z(L_{c+1} g(z))' + (c+1) L_{c+1} g(z)} \\
& & = Re \fr{ \fr{z(L_{c+1} zf'(z))'}{L_{c+1} g(z)} + (c+1) \fr{z(L_{c+1} f(z))'}{L_{c+1} g(z)}}{\fr{z(L_{c+1} g(z))'}{L_{c+1} g(z)} + c + 1} .
\eenrr
In view of (2.13) and considering Theorem 1 we have $g \in S^*_{c+1} (\la)$ and $\fr{z(L_{c+1} g(z))'}{L_{c+1} g(z)} = (1 - \la) Q(z) + \la$ where $Re(Q(z)) > 0, z \in U$, also according to (2.16) we have
$$L_{c+1} zf'(z) = L_{c+1} g(z) [(1 - \b) p(z) + \b] . \eqno{(2.20)}$$
Differentiating logarithmically with respect to $z$ from both sides of (2.20) gives
$$ \fr{z(L_{c+1} zf'(z))'}{L_{c+1} g(z)} = (1 - \b) zp'(z) + [(1-\la) Q(z) + \la] [(1 - \b) p(z) + \b].\eqno{(2.21)} $$
However,
\benrr
& & Re \left\{ \fr{z(L_cf(z))'}{L_cg(z)} \right\}\\& &  < Re \fr{(1 - \b)zp'(z) + [(1-\la)Q(z)+\la] [(1-\b)p(z)+\b] + (c+1) [(1-\b)p(z)+\b]}{(1 - \la) Q(z) + \la + c + 1} \\
& & = Re \{ (1 - \b) p(z) + \b\} + Re \fr{(1 - \b) zp'(z)}{(1 - \la) Q(z) + \la + c + 1} .
\eenrr
Equivalently
$$ Re \left\{ \fr{z(L_cf(z))'}{L_cg(z)} - \b \right\} < Re \left\{ (1 - \b) p(z) + \fr{(1 - \b) zp'(z)}{(1 - \la)Q(z) + \la + c + 1} \right\} .\eqno{(2.22)}$$
By considering the function $\psi(u, v)$ as
$$\psi(u, v) = (1 - \b)u + \fr{(1 - \b)v}{(1 - \la) Q(z) + \la + c + 1}$$
and noting that $Re(Q(z)) > 0$ we can easily verify that the function $\psi$ is a continuous function in $D = \mathbb{C} \times  \mathbb{C}$ and $Re\{ \psi(1, 0)\} > 0,$ also if $v_1 \le - \fr{1}{2}  (1 + u^2_2)$ then we have
\benrr
& & Re \{ \psi(iu_2, v_1) = Re \left\{(1 - \b) iu_2 + \fr{(1-\b)v_1}{(1-\la) Q(z) + \la + c + 1} \right\} \\
& & = Re \left\{ \fr{(1-\b)v_1[\la +c+1 + (1-\la) Re(Q(z)) - i(1-\la) I_m(Q(z))]}{[\la + c + 1 + (1-\la) Re(Q(z))]^2 + [(1-\la) I_m(Q(z))]^2} \right\} \\
& & = \fr{(1 - \b)v_1[\la + c + 1 + (1 - \la) Re(Q(z))]}{[\la + c + 1 + (1-\la) Re(Q(z))]^2 + [(1-\la) I_m(Q(z)]^2} \\
& & \le \fr{-(1-\b) (1 + u^2_2) [\la + c + 1 + (1-\la)
Re(Q(z))]}{[\la + c + 1 + (1-\la) Re(Q(z))]^2 +
[(1-\la)I_m(Q(z))]^2} < 0. \eenrr Finally since in view of (2.22)
we have $Re\{ \psi(p(z), zp'(z))\} > 0$ therefore Lemma 1.1 gives
$Re(p(z)) > 0, z \in U$ and the proof is complete.

The proof of part (ii) is similar to part (i) and we omit it.

By the same method used in Theorem 6 we can prove the next theorem.

\noi {\bf Theorem 2.9} : (i) If for every $f \in A$ and $g \in
C_c(\la)$ we have
$$Re \left\{ \fr{z \fr{d}{dz} \left( \fr{(L_c zf'(z))'}{(L_cg(z))'} \right)}{\fr{z(L_cg(z))''}{(L_cg(z))'} + c + 1 } \right\} > 0 $$
and
$$Re \left\{ \fr{zg'(z)}{g(z)} - \fr{z(L_cg(z))'}{L_cg(z)} \right\} > 0 $$
then $K^*_c(\b, \la) \subset K^*_{c+1} (\b, \la)$.

(ii) If for every $f \in A$ and $g \in C_{c+1}(\la)$ we have
$$Re \left\{ \fr{z \fr{d}{dz} \left( \fr{(L_{c+1} zf'(z))'}{(L_{c+1}g(z))'} \right)}{\fr{z(L_{c+1}g(z))''}{(L_{c+1}g(z))'} + c + 1 } \right\} > 0 $$
and
$$Re \left\{ \fr{zg'(z)}{g(z)} - \fr{z(L_{c+1}g(z))'}{L_{c+1}g(z)} \right\} > 0 $$
then $K^*_{c+1}(\b, \la) \subset K^*_c (\b, \la)$.

\noi {\bf Theorem 2.10} : If $- \la \le c \le 1 - 2\la$ then $f
\in S^*_\si(\la)$ implies $I^\si f \in S^*_c(\la)$.

\noi {\bf Proof} : Suppose that $f \in S^*_\si(\la)$ and set
$$ \fr{z(L_cI^\si f(z))'}{L_c I^\si f(z)} = \fr{1 + (1 - 2\la) w(z)}{1 - w(z)}, ~~ z \in U \eqno{(2.23)}$$
where $w(z)$ is analytic or meromorphic in $U$ with $w(0) = 0$.  By using (1.8) and (2.23) we obtain
$$ \fr{I^\si f(z)}{L_c I^\si f(z)} = \fr{c+1 + (1-c-2\la) w(z)}{(c+1) (1 - w(z))} .\eqno{(2.24)}$$
Differentiating logarithmically both sides of (2.24) with respect to $z$ gives
$$ \fr{z(I^\si f(z))'}{I^\si f(z)} = \fr{1 + (1 - 2\la) w(z) + zw'(z)}{1 - w(z)} + \fr{(1-c - 2\la) zw'(z)}{c+1 + (1 - c - 2\la) w(z)} $$
Now we assert that $|w(z)| < 1, z \in U$, if not then there exists a point $z_0 \in U$ such that $\max\limits_{|z| \le |z_0|} |w(z)| = |w(z_0)| = 1$ therefore by Jacks' Lemma we have $z_0 w'(z_0) = kw(z_0), k \ge 1$.  So we have
\benrr
& & Re \left\{ \fr{z_0(I^\si f(z_0))'}{I^\si f(z_0)} - \la \right\} \\
& & = Re \left\{ \fr{1 + (1 - 2\la + k)e^{i\th}}{1 - e^{i\th}} + \fr{(1 - c - 2\la) ke^{i\th}}{c + 1 + (1-c-2\la)e^{i\th}} - \la \right\} \\
& & = \fr{-2k(1 - \la) (c + \la)}{(1+c)^2 + 2(1+c) (1-c-2\la) \cos \th + (1 - c - 2\la)^2} \le \fr{-k(c+\la)}{2(1-\la)} \le 0.
\eenrr
This contradicts our hypothesis $f \in S^*_\si(\la)$ thus $|w(z)| < 1, z \in U$ and by cosidering (2.23) we conclude that $I^\si f \in S^*_c(\la)$.

\noi {\bf Corollary 2.11} : If $- \la < c < 1 - 2\la$ and $f \in
C_\si(\la)$ then $I^\si f \in C_c (\la)$.

\noi {\bf Proof} : We have \\$f \in C_\si(\la) \Leftrightarrow zf' \in S^*_\si(\la) \La I^\si(zf') \in S_c^*(\la) \Leftrightarrow z(I^\si f)' \in S^*_c(\la) \Leftrightarrow I^\si f \in C_c(\la)$.

\noi {\bf Theorem 2.12} : Let $- \la \le c, 0 \le \la < 1$.  If $f
\in A$ and $\fr{z(L_cI^\si f(z))'}{L_c I^\si f(z)} \ne \la, z \in
U$ then $f \in ST_\si(\eta, \la)$ implies that $I^\si f \in
ST_c(\eta, \la)$.

\noi {\bf Proof} : Let $f \in ST_\si (\eta, \la)$ and put
$$ \fr{z(L_c I^\si f(z))'}{L_c I^\si f(z)} = \la + (1 - \la) p(z) \eqno{(2.25)}$$
where $p(z) = 1 + \sum\limits^\iny_{n=1} c_nz^n$ and $p(z) \ne 0, z \in U$.  By considering (1.8) and (2.25) we have
$$ (c+1) \fr{I^\si f(z)}{L_c I^\si f(z)} = c + \la + (1 - \la) p(z) \eqno{(2.26)}$$
Differentiating logarithmically with respect to $z$ from both sides of (2.26) gives
$$ \fr{z(I^\si f(z))'}{I^\si f(z)} - \la = (1 - \la) p(z) + \fr{(1 - \la) zp'(z)}{c + \la + (1 - \la)p(z)} .$$
Suppose that there eixsts a point $z_0 \in U$ such that
$|arg(p(z))| < \fr{\pi}{2} \eta (|z| < |z_0|)$ and $|arg(p(z_0))|
= \fr{\pi}{2}\eta$ then by Lemma 1.2 we have
$\fr{z_0p'(z_0)}{p(z_0)} = ik\eta$ and $p(z_0)^{1/\eta} = \pm~~
ir (r > 0)$.  If $p(z_0)^{1/\eta} = ir$ then \benrr
& & \fr{z_0 (I^\si f(z_0))'}{I^\si f(z_0)} - \la = (1 - \la) p(z_0) \left[ 1 + \fr{\fr{z_0p'(z_0)}{p(z_0)}}{c + \la + (1 - \la) p(z_0)} \right] \\
& & = (1 - \la) r^\eta e^{i \fr{\pi}{2} \eta} \left[ 1 + \fr{ik\eta}{c + \la + (1 - \la)r^\eta e^{i \fr{\pi}{2} \eta}} \right] \\
& & = \fr{\pi}{2} \eta + arg \left\{ 1 + \fr{ik\eta}{c + \la + (1 - \la)r^\eta e^{i \fr{\pi}{2} \eta}} \right\} \\
& & = \fr{\pi}{2} \eta \\
& & + \tan^{-1} \left\{ \fr{k\eta [c + \la + (1 - \la)r^\eta \cos \fr{\pi}{2} \eta]}{(c + \la)^2 + 2(c + \la)(1 - \la)r^\eta \cos \fr{\pi}{2} \eta + (1-\la)^2 r^{2\eta} + k\eta (1 - \la)r^\eta \sin \fr{\pi}{2} \eta} \right\} \\
& & \ge \fr{\pi}{2} \eta ~~~ (\mbox{Because} ~~~ k \ge \fr{1}{2} (r + \fr{1}{r}) \ge 1)
\eenrr
which contradicts our assumption $f \in ST_\si(\eta, \la)$.  By the same method we get a contradiction for the case $p(z_0)^{1/\eta} = -ir (r > 0)$, therefore we have $|arg(p(z))| < \fr{\pi}{2} \eta, z \in U$ and in view of (2.14) we conclude that $I^\si f \in ST_c(\eta, \la)$.

\noi {\bf Corollary 2.13} : Let $c \ge \la, 0 \le \la < 1$.  If $f
\in A$ and $\fr{(z(L_c I^\si f(z))')'}{(L_c I^\si f(z))'} \ne
\la, z \in U$ then $f \in CV_\si (\eta, \la)$ implies that $I^\si
f \in CV_c(\eta, \la)$.

We claim the similar results may be hold for meromorphic $p$-valent
functions with alternating coefficient. For more information see
[2].

\newpage

\bc
{\bf References} \\
\ec

\begin{verse}
[1] S. D. Bernardi, Convex and starlike univalent functions, Trans. Amer. Math. Soc., 135 (1969), 429-446.

[2] A. Ebadian, S. Shams and Sh. Najafzadeh,  Certain inequalities
for $p$-valent meromorphic functions with alternating coefficients
based on integral operator. Aust. J. Math. Anal. Appl. 5 (2008), no.
1, Art. 10, 5 pp.

[3] I. B. Jung, Y. C. Kim and H. M. Srivastava, The Hardy space of
analytic functions associated with certain one-parameter families of
integral operators, J. Math. Anal. Appl., 176 (1993), 138-147.

[4] J. L. Li, Some properties of two integral operators, Soochow. J.
Math., 25 (1999), 91-96.

[5] R. J. Libera, Some classes of regular functions, Proc. Amer.
Math. Soc., 16 (1965), 755-758.

[6] J. L. Liu, A linear operator and strongly starlike functions, J.
Math. Soc. Japan, Vol. 54, No. 4 (2002), 975-981.

[7] ---------- , Some applications of certain integral operator,
Kyungpook Math. J. 43(2003), 21-219.

[8] A. E. Livingston, On the radius of univalence of certain
analytic functions, Proc. Amer. Math. Soc., 17(1996)352-357.

[9] S. S. Miller and P. T. Mocanu, Second order differential
inequalities in the complex plane, J. Math. Anal. Appl., 65 (1978),
289-305.

[10] K. I. Noor, On quasi-convex functions and related topics,
Internat. J. Math. Math. Sci., 10 (1987), 241-258.

[11] M. Nunokawa, S. Owa, H. Saitoh, A. Ikeda and N. Koike, Some
results for strongly starlike functions, J. Math. Anal. Appl., 212
(1997), 98-106.

[12] B. A. Uralegaddi and C. Somanatha, Certain integral operators
for starlike functions, J. Math. Res. Expo., 15 (1995), 14-16.
\end{verse}

\end{document}